\documentclass[psfig,epsfig]{article}
\begin{document}

\newcommand{\ea}{\mbox{{\bf a}}}
\newcommand{\eu}{\mbox{{\bf u}}}
\newcommand{\ueu}{\underline{\eu}}
\newcommand{\ueo}{\overline{u}}
\newcommand{\oeu}{\overline{\eu}}
\newcommand{\ew}{\mbox{{\bf w}}}
\newcommand{\ef}{\mbox{{\bf f}}}
\newcommand{\eF}{\mbox{{\bf F}}}
\newcommand{\eC}{\mbox{{\bf C}}}
\newcommand{\en}{\mbox{{\bf n}}}
\newcommand{\eT}{\mbox{{\bf T}}}
\newcommand{\eL}{\mbox{{\bf L}}}
\newcommand{\eV}{\mbox{{\bf V}}}
\newcommand{\eU}{\mbox{{\bf U}}}
\newcommand{\ev}{\mbox{{\bf v}}}
\newcommand{\eve}{\mbox{{\bf e}}}
\newcommand{\uev}{\underline{\ev}}
\newcommand{\eY}{\mbox{{\bf Y}}}
\newcommand{\eK}{\mbox{{\bf K}}}
\newcommand{\eP}{\mbox{{\bf P}}}
\newcommand{\eS}{\mbox{{\bf S}}}
\newcommand{\eJ}{\mbox{{\bf J}}}
\newcommand{\eB}{\mbox{{\bf B}}}
\newcommand{\leb}{{\cal L}^{n}}
\newcommand{\eI}{{\cal I}}
\newcommand{\eE}{{\cal E}}
\newcommand{\hen}{{\cal H}^{n-1}}
\newcommand{\eBV}{\mbox{{\bf BV}}}
\newcommand{\eA}{\mbox{{\bf A}}}
\newcommand{\eSBV}{\mbox{{\bf SBV}}}
\newcommand{\eBD}{\mbox{{\bf BD}}}
\newcommand{\eSBD}{\mbox{{\bf SBD}}}
\newcommand{\ecs}{\mbox{{\bf X}}}
\newcommand{\eg}{\mbox{{\bf g}}}
\newcommand{\paromega}{\partial \Omega}
\newcommand{\gau}{\Gamma_{u}}
\newcommand{\gaf}{\Gamma_{f}}
\newcommand{\sig}{{\bf \sigma}}
\newcommand{\gac}{\Gamma_{\mbox{{\bf c}}}}
\newcommand{\deu}{\dot{\eu}}
\newcommand{\dueu}{\underline{\deu}}
\newcommand{\dev}{\dot{\ev}}
\newcommand{\duev}{\underline{\dev}}
\newcommand{\weak}{\rightharpoonup}
\newcommand{\weakdown}{\rightharpoondown}
\renewcommand{\contentsname}{ }

\title{Energy concentration and brittle crack propagation}
\date{ }


\maketitle

\begin{abstract}

The purpose of this paper is to fill the gap between the 
classical treatment of brittle fracture mechanics and the new idea of considering the 
crack evolution as a free discontinuity problem. Griffith and Irwin 
criterions of crack propagation are studied and transformed in order to be no longer dependent 
on any prescription of the geometry of the crack during its evolution. 
The inequality contained in theorem 6.1. represents the link between generalized Irwin and  
Griffith criterions of brittle crack propagation. The physical meaning of this inequality 
is explained in the last section. 

\end{abstract}

\tableofcontents

\newpage

\section*{Note to the reader}

This is a paper which has been submitted for publication in a well-known journal in 1997. The referee 
 did not want to accept  the paper  unless substantial modifications are made. In my opinion the suggested  modifications (excepting ortographic,  typographical or minor mathematical ones) were 
 against the philosophy of the paper, which is: brittle fracture propagation is a geometrical evolution 
 problem, therefore geometrical treatment is highly significant, both from mechanical and mathematical point of view. 
 
 As this (eventually unpublished) paper circulated in manuscript, I think it is safe to depose it in the arXiv, for further reference and comparison with other works in this narrow domain. 
 
 No modifications were made to the original file, excepting some preprint references in the bibliography, which were published since.

\section{Introduction}
\indent

Brittle fracture mechanics studies the evolution of the cracks in elastic bodies. 
Since the medium under consideration is supposed to remain elastic, the problem of 
brittle crack propagation concerns the evolution of the geometric support of the 
crack. We are facing here a problem placed in between geometry and mechanics. 

The field of research in brittle fracture mechanics has appeared in 1920 with  
Griffith's paper [G]. 
Among the basic references in this field we find:  Irwin [I], Eshelby [Es], 
Gurtin [Gu1], [Gu2]. 

Typically for classical fracture mechanics is that  
the geometrical nature of brittle crack evolution problem is obscured 
 by the assumption of prescribed geometry of the crack. There are very few 
papers which do not make this assumption; as an example we cite here Ohtsuka [Oht1---4] and 
Stumpf \& Le [StLe].   

  In the last decade a new current of ideas has emerged, starting with the article of 
Mumford \& Shah [MS], submitted for publication in 1986. It is the first time, to our knowledge, 
when the crack itself is considered as the unknown of the problem of brittle crack appearance. 
Some new mathematical results of De Giorgi and Ambrosio ([DGA], [A1], [A2]) in the field 
of geometric measure theory, set the energetic model of brittle crack appearance proposed by 
Mumford and Shah in the functional space $\eSBV$, of special functions with bounded variation. 

The strange unknown of  the problem --- the crack --- is replaced by a more familiar one: 
the displacement field of the cracked body. This field is allowed to be discontinuous, being a 
special function with bounded variation.  Recent papers study the space $\eSBD$ of special 
functions with bounded deformation, as Bellettini, Coscia \& Dal Maso [BCDM] and Ambrosio, 
Coscia \& Dal Maso [ACDM], naturally 
associated with the expression of free energy potential of an elastic body 
suffering small deformations. 

It is now clear that a promising point of view concerning crack evolution is to consider it 
as a free discontinuity problem. The purpose of this paper  is to fill the gap between the 
classical treatment of brittle fracture mechanics and these new ideas. We present here 
some of the results obtained in [Bu].  The Griffith and Irwin 
criterions of crack propagation are studied and transformed in order to be no longer dependent 
on any prescription of the geometry of the crack during its evolution. 

The content of the paper is described further. In the section "Preliminaries and notations"  
some basic constitutive assumptions are made. The following section, "Statics of a fractured 
body", contains a mathematical treatment of the equilibrium of an elastic cracked body in a 
functional setting compatible with the space $\eSBV$. This section contains also a brief
description of the Dirichlet-to-Neumann map of a cracked body. For more information about this 
notion we send the reader to  [C], [SU1,2]. 

The section "Kinematics of crack propagation"  starts with the introduction of a set  
of  reasonably smooth endomorphisms of $\Omega$, the reference configuration of the body, 
useful in the sequel. The notion of smooth crack propagation is 
introduced, in various forms (definitions  4.1. to  4.3.). A smooth crack propagation can be seen 
as a curve of endomorphisms of $\Omega$, namely $t \mapsto \phi_{t}$. The initial crack $K$, which
is a surface with boundary in $\Omega$, becomes at the moment $t$ the actual crack 
$K_{t} = \phi_{t}(K)$. We mention that in [StLe] variations of the crack by diffeomorphisms are 
considered; also, in [Oht1---4], a curve of diffeomorphisms $t \mapsto \phi_{t}$ 
is associated to a curve of cracks $t \mapsto K_{t}$. 

In the section "Crack propagation criterions"  Griffith and Irwin criterions of brittle crack 
propagation are reformulated for smooth crack propagation. This is accomplished by proposition 5.1., 
 definitions 5.1. (generalized Rice's J integral) and 5.2 (concentration coefficients of the
elastic energy).  

The section "Energy concentration and crack propagation" begins with the introduction of the
measure $\mid K2 \mid(\eu)$ (definition 6.1.) associated to a field $\eu$ of displacements. 
This measure describes the distribution in the body of the energy release rate due to 
the crack propagation. In definition 6.2. is introduced, for the same arbitrary field of 
displacements $\eu$,  the upper energy concentration coefficient
as a measure, namely $CM^{+}(\eu)$. The main result of the paper is theorem 6.1., which 
describes the relations between these two measures. 

A byproduct  of the proof of theorem 6.1., 
with important physical significance,  is presented in the last section.

\section{Preliminaries and notations}
\indent

Let $\Omega \subset R^{n}$ represent the reference configuration of an elastic body. 
 The space dimension $n$ equals 2 or 3 and 
$\Omega$ is an open bounded set with piecewise smooth boundary $\paromega$. 
The volume density of the free energy is a function $w = w (\nabla \eu)$. $\eu$ denotes 
the displacement of the body from the reference configuration. We suppose that the function $w$
 depends only on $\epsilon(\eu)$, the symmetric part 
of $\nabla \eu$. In this paper we work with a quadratic expression of $w$:
\begin{equation}
w(\nabla \eu) =  \frac{1}{2} \eC \nabla \eu : \nabla \eu   \ \ \ \ \  .
\label{ipcon}
\end{equation}
The tensor $\eC$ has the symmetries: 
\begin{equation}
C_{ijkl} \ = \ C_{jikl} \ = \ C_{klij} \ \ \ \ \ \ .
\label{hyp}
\end{equation}
Most of the facts herein are true whenever $w$ is a convex, $C^{2}$ function, satisfying 
the following growth condition: there are two positive constants $c,C$ such that 
\begin{equation}
 \forall E \in R^{n \times n} \ , \ c \ \mid  E^{sym} \mid^{2}  \ \leq \ w(E) 
\ \leq C \mid E \mid^{2} \ \ \ \ \ \left( E^{sym} := \frac{1}{2} ( E + E^{T}) \right) \ \ .
\label{grow}
\end{equation} 
We will restrict however our attention to the case when the energy function $w$  
satisfies (\ref{grow}), has the form 
(\ref{ipcon}) and   $\eC$ satisfies (\ref{hyp}). 
  
The first Piola-Kirchoff stress tensor is defined by:
$$\sigma_{ij} = \sigma_{ij} (\nabla \eu) =  \frac{d w}{d A_{ij}}(\nabla \eu) \ \ \ \ \ .$$
The body evolves in a quasi-static manner, in the absence of volume forces, so at any moment 
the stress tensor $\sigma$ is divergence free.

 The 
Hausdorff $k$-dimensional measure is denoted by ${\cal H}^{k}$. ${\cal H}^{2}$ is the area 
measure and ${\cal H}^{1}$ is the length measure. 
In order to simplify the denominations we shall often call $\hen$ the area measure. 
The volume measure is 
${\cal H}^{n} =  \leb$, the Lebesgue measure. 

Any crack in the
body is seen as a "crack set" (see Ball [Ba]): a closed countably rectifiable set 
$K \subset \overline{\Omega}$ with finite area.  
 Any smooth hyper-surface with bounded mean curvature and with boundary in $\overline{\Omega}$ is 
called a "smooth crack set". A smooth crack set it is therefore a crack set which 
carries the geometric structure of a manifold with boundary. The boundary of the smooth crack set $K$, 
denoted by $\partial K$, represents the edge of the
smooth crack set. The edge $\partial K$ is a $n-2$ piecewise smooth surface without boundary.

In this paper we work with smooth crack sets $K$ with the property:
$${\cal H}^{n-1} \left( K \cap \partial \Omega \right) \ = \ 0 \ \ .$$
This assumption has the following meaning: the intersection of the smooth crack set with the 
exterior boundary of the body $\partial \Omega$ is at most a reunion of curves. 

The space $\eSBV(\Omega,R^{n})$ of special functions with bounded variation was introduced by 
De Giorgi and Ambrosio in the study of a class of free discontinuity problems ([DGA], [A1], [A2]). 
For any function $\eu \in L^{1}(\Omega,R^{n})$ let us denote by $D\eu$ the distributional 
derivative of $\eu$ seen as a vector measure. The variation of $D\eu$ is a scalar measure 
defined like this:  
for any Borel measurable subset $B$ of $\Omega$ the variation of $D\eu$ over $B$ is 
\begin{displaymath}
\mid D\eu \mid  (B) \ = \  sup \ \left\{ \sum^{\infty}_{i=1} \mid D\eu  (A_{i}) \mid 
\mbox{ : } \cup_{i=1}^{\infty} A_{i} \subset B \ , \ A_{i} \cap A_{j} = \emptyset \ \ \forall 
i \not = j \right\} \ \ .
\end{displaymath}
The Lebesgue set of $\eu$ is the set of points where $\eu$ has approximate limit. The complementary 
set is a $\leb$ negligible set denoted by $\eS_{\eu}$.  

 The space $\eSBV(\Omega,R^{n})$ is defined as follows: 
$$\eSBV(\Omega,R^{n}) \ = \ \left\{ \eu \in L^{1}(\Omega,R^{n}) \mbox{ : } \mid D\eu \mid (\Omega)
 < + \infty 
\ , \ \mid D^{s}\eu \mid (\Omega \setminus \eS_{\eu}) = 0 \right\} \ .$$
Let us define the following Sobolev space  associated to the crack set $K$ (see [ABF]):
$$W^{1,2}_{K} \ = \ \left\{ \eu \in \eSBV(\Omega,R^{n}) \mbox{ : } \int_{\Omega} 
\mid \nabla \eu \mid^{2} 
\mbox{ d}x + \int_{K} [\eu]^{2} \mbox{ d} \hen  <  + \infty \ , \ \mid D^{s}\eu \mid \ll \hen_{|_{K}} 
\right\} \ .$$
It has been proved in [DGCL] the following equality:
$$W^{1,2}_{K}(\Omega,R^{n}) \cap L^{\infty}(\Omega,R^{n}) = W^{1,2}(\Omega \setminus K, R^{n}) 
\cap L^{\infty}(\Omega,R^{n}) \ \ \ .$$

We shall need further a norm with physical dimension on the Sobolev space $$W^{p,2,\infty}(\Omega,R^{n}) = W^{p,2}(\Omega,R^{n}) \cap 
L^{\infty}(\Omega,R^{n}) \ \ \ \ .$$ Consider therefore $\lambda$ as an unit of length. 
We take the following  norm on $W^{p,2}(\Omega,R^{n})$~:  
$$ \| \eu \|_{p,2} = \| \eu \|_{L^{2}} \ + \ \lambda \ \| \nabla \eu \|_{L^{2}} \ + \  ... \ + \ 
\lambda^{p} \ \| \nabla^{p} \eu \|_{L^{2}} \ \ \ \ .$$
This norm induces on $W^{p,2,\infty}(\Omega,R^{n})$ the norm  
$$\| \eu \|_{p,2,\infty} = \ max ( vol(\Omega) \ \| \eu \|_{\infty}, \| \eu \|_{p,2}) \ \ \ ,$$ 
where $vol(\Omega)$ is the volume of $\Omega$ measured with the unit $\lambda^{n}$.

\section{Statics of a fractured body}
\indent

\subsection{Functional assumptions}
\indent

Let us consider an elastic body $\Omega$ with a smooth crack set $K$. The displacement 
$\eu^{0} \in H^{\frac{1}{2}}(\paromega,R^{n}) \cap L^{\infty}(\paromega,R^{n})$ is imposed 
on the exterior boundary $\paromega$. 

Therefore the equilibrium displacement $\eu$ is a special function with bounded variation. 
>From the Calderon \& Zygmund [CZ] decomposition theorem we obtain  the following expression of 
$D\eu$, the distributional derivative of $\eu$ seen as a measure: 
$$D\eu \ = \  \nabla \eu (x) \mbox{ d}x \ + \  [\eu ] \otimes \en \mbox{ d}\hen_{|_{K}} \ \ \ .$$ 
We deduce from here the Stokes formula:
\begin{equation}
\langle D\eu , \phi \rangle = \int_{\Omega} \nabla \eu \cdot \phi \mbox{ d}x + 
\int_{K} [\eu] \otimes \en \cdot \phi \mbox{ d}\hen \mbox{ \ \ } \forall \phi 
\in C^{\infty}_{0}(\Omega , R^{n \times n}) \ \ .
\label{stokes}
\end{equation}
 This formula is the mathematical 
expression of the fact that are no concentrated forces on the edge of the crack. 

Any element of $W^{1,2}_{K}$, compatible with the boundary condition, is called an admissible 
displacement. The family of admissible stresses will be defined in the sequel.

The jump of a field $\sigma$ across $K$ is $[\sigma] = \sigma^{+} - \sigma^{-}$. Del Piero
considers 
in [DP] the following set of admissible stresses for a fractured media: 
$$Y (\Omega) = \left\{ \sigma = \sigma^{T} \in L^{2}(\Omega,R^{n \times n}) \mbox{ : } div \sigma  
= 0 \right\} \ \ .$$ 
He proves that if $\sigma \in Y(\Omega)$ then $ \sigma^{+}  \en = \sigma^{-} \en$. This fact 
allows us to use the notation $\sigma \en$ without confusion.   

We connect further the definition of an admissible displacement field with the definition of an 
admissible stress field. 
 
\vspace{.5cm}

{\bf Definition 3.1.1.} {\it Let $\eu^{0} \in H^{\frac{1}{2}}(\paromega,R^{n}) \cap L^{\infty}(\paromega,R^{n})$ 
represent an imposed displacement on the exterior boundary of the body $\Omega$ and let $K$ be a 
crack set in $\Omega$. The set of admissible displacements with respect to $\eu^{0}$ and $K$ is 
\begin{equation}
W^{1,2,\infty}_{K}(\eu^{0}) \ = \ \left\{ \eu \in W^{1,2}(\Omega \setminus K, R^{n}) 
\cap L^{\infty}(\Omega,R^{n}) \mbox{ : } \eu = \eu^{0} \mbox{ on } \paromega \ \right\}
\label{disadm}
\end{equation} 
and the set of admissible stresses with respect to  $K$ is} 
\begin{equation}
Y_{K} (\Omega) = \left\{ \sigma = \sigma^{T} \in L^{2}(\Omega,R^{n \times n}) \mbox{ : } div \sigma  
= 0 \ , \ \sigma \en = 0 \mbox{ on } K \right\} \ \ .
\label{stadm}
\end{equation}

\vspace{.5cm}

On the crack set the body is force free, i.e. 
$\sigma^{+,-} \en = 0$ on $K$, where $\en$ is 
the field of normals on $K$ and $\sigma^{+}$, $\sigma^{-}$ are the traces of $\sigma$ on the sides 
of $K$. The sign convention is taken such that the  field of normals  to $K$, $\en$, 
points to the "$+$" side of $K$. There are no concentrated forces 
on the edge of the crack.

\subsection{Dirichlet-to Neumann map}
\indent

The equilibrium displacement of the body, when the boundary displacement $\eu^{0}$ is imposed, is 
solution of the problem:
\begin{equation}
\left\{ \begin{array}{ll}
 div \ \sigma(\nabla \eu) \ = \ 0 & \mbox{ in } \Omega \setminus K \\
 \sigma(\nabla \eu)^{+} \en \ = \ \sigma(\nabla \eu)^{-} \en \ = \ 0 & \mbox{ on } K \\
\eu \ = \ \eu^{0} & \mbox{ on } \paromega  \ \ \ \ . 
\end{array} \right.
\label{chilu} 
\end{equation}

 The equilibrium displacement 
$\eu$ minimizes the functional 
$$E(\ev) \ = \ \int_{\Omega} w(\nabla \ev) \mbox{ d}x$$
defined over all $\ev \in W^{1,2}(\Omega \setminus K, R^{n})$, $\ev = \eu^{0}$ on $\paromega$. 
The solution of the minimization problem is unique to an arbitrary piecewise rigid displacement 
equal to 0 on $\paromega$. Also, $\eu$ is  smooth and essentially bounded.

$Y_{K} (\Omega)$ is a subset of $Y(\Omega)$; in the definition the quantity $\sigma \en$ does not 
depend on the choice of $\en$.

The polar of the convex functional $W(\epsilon) = \int_{\Omega} w(\epsilon) \mbox{ d}x$ is (see 
Moreau [M]): 
 $$W^{*}(\sigma) = \sup \left\{ \langle \sigma , \epsilon \rangle - W(\epsilon) \mbox{ : } 
\epsilon \in L^{2}(\Omega,R^{n \times n}_{sym}) \right\} \ \ .$$

For any admissible stress $\sigma \in Y_{K}(\Omega)$ and for any
minimizer $\eu$ of the functional $E$ over the class of admissible displacements, the following 
inequality is true: 
\begin{equation}
\int_{\paromega} (\sigma \en) \cdot \eu_{0} \mbox{ d}\hen - \ W^{*}(\sigma) \leq E(\eu) \ \ .
\label{moreau}
\end{equation}

The inequality becomes an equality if and only if $\sigma = \sigma (\nabla \eu)$.

\vspace{2.cm}

For a given crack set $K$ the Dirichlet-to-Neumann (or the response) map associated to the domain 
$\Omega$ with the crack $K$ inside can be defined. This map is defined in the 
following way: let  $\eu^{0}$ be given 
on the boundary of $\Omega$. The solution of the Dirichlet problem 
\begin{equation}
\left\{ \begin{array}{ll}
 div \ \sigma(\nabla \eu ) \ = \ 0 & \mbox{ in } \Omega \setminus K \\
 \eu \ = \ \eu^{0} & \mbox{ on } \paromega \\
  \sigma(\nabla \eu)  \en \ = \ 0 & \mbox{ on (both sides of)  } K
\end{array} \right. 
\label{diric}
\end{equation} 
will be denoted by $\eu(K, \eu^{0})$. This solution is uniquely determined to a rigid 
displacement equal to 0 on $\paromega$ hence $\sigma(K,\eu^{0})$ is unique. 

It is obvious that 
$\eu(K, \eu^{0})$ is also solution of the Neumann problem 
\begin{equation}
\left\{ \begin{array}{ll}
 div \ \sigma(\nabla \eu ) \ = \ 0 & \mbox{ in } \Omega \setminus K \\
 \sigma(\nabla \eu)  \en \ = \ \sigma(\nabla \eu(K,\eu^{0})) \en & \mbox{ on } \paromega \\
  \sigma(\nabla \eu) \en \ = \ 0 & \mbox{ on (both sides of)  } K  \ \ .
\end{array} \right. 
\label{neuma}
\end{equation}
Therefore the function  $\eu^{0} \mapsto \sigma(\nabla \eu(K,\eu^{0})) \cdot \en$ maps naturally 
a Dirichlet boundary condition to a Neumann boundary condition. The Dirichlet-to-Neumann map is 
the linear application: 
$$\eT(K) : H^{\frac{1}{2}}(\paromega) \rightarrow H^{-\frac{1}{2}}(\paromega) \ \ ,$$ 
\begin{equation}
\langle \eT(K) \eu^{0} , \ev^{0} \rangle \ = \ 
\int_{\paromega} \sigma(\nabla \eu(K,\eu^{0})) \en \cdot  \ev \mbox{ d}\hen  \ \ .
\label{dndef}
\end{equation}
Here $\langle \ , \ \rangle$ means the duality product.The Dirichlet-to-Neumann map is 
continuous due to the well known continuous dependence 
of the solutions of the problems (\ref{diric}) and (\ref{neuma}) with respect to boundary data. 
The symmetries of the elasticity tensor $\eC$ make this map self adjoint.

If $\eu^{0}$ is essentially bounded then 
the solution of the problem (unique to a rigid displacement equal to 0 on $\paromega$) 

\hspace{2.cm} {\it minimize 
$ \int_{\Omega} w(\nabla \eu) \mbox{ d}x $
 over all $\eu \in W^{1,2,\infty}_{K} (\eu^{0})$} 
  
\noindent
is $\eu(K, \eu^{0})$, hence 
$$\frac{1}{2} \ \langle \eT(K) \eu^{0} , \eu^{0} \rangle \ = \ min \ 
\left\{  \int_{\Omega} w(\nabla \eu) \mbox{ d}x \mbox{ : }  \eu \in 
W^{1,2,\infty}_{K} (\eu^{0}) \right\} \ \ .$$ 





\section{Kinematics of crack propagation}
\indent

Let us consider a hyper-elastic body with the reference configuration $\Omega$. At the moment $t=0$ in the body 
there is the crack set $K$. The evolution of the crack is a curve $t \mapsto K_{t}$ such that  $K_{0} = K$ and  
for any $t < t'$ $K_{t} \subset K_{t'}$.

We are searching for an easy way to describe the evolution of the crack set. If the dimension 
of the space is $n=2$ then an obvious way to describe the 
evolution of the crack set is to consider only the evolution of the  edge of $K$,i.e. $\partial K$,  which is formed
by a finite number of points. 
In the general case $n=3$ this corresponds to the evolution $t \mapsto \partial K_{t}$ where 
$\partial K_{t}$ is  a curve. Because the evolution of a surface in the 
three-dimensional space is conceptually simpler than the evolution of a curve, the before 
mentioned choice of modeling the crack propagation  does not simplify our problem. We prefer instead to think at 
the crack set $K_{t}$ as a deformation of the initial crack set $K$ by a diffeomorphism 
$\phi_{t} : \Omega \rightarrow \Omega$. We replace therefore the curve $t \mapsto K_{t}$ with 
the curve $t \mapsto \phi_{t}$, satisfying the condition $\phi_{t}(K) = K_{t}$. (However, not 
any curve $t \mapsto K_{t}$ has an associated curve $t \mapsto \phi_{t}$. As a counterexample 
we may think at cracks 
which after a time develop new branches.) 

The curve $t \mapsto \phi_{t}$ lies in a set of diffeomorphisms.

Let us consider the following family of diffeomorphisms:
\begin{equation}
M_{c} \ = \  \left\{  
\phi \in C^{\infty}(\Omega,\Omega) \cap 
W^{s,2,\infty}(\Omega,R^{n}) \mbox{  :  } \right. 
\label{mec}
\end{equation}
$$\left. \phi \mbox{ is a diffeomorphism  and } supp \ (\phi - 1_{\Omega} ) 
\subset \Omega \ \right\}$$

The set $M_{c}$ can be conveniently completed to a topological group, which has 
a geometrical structure.

We consider, inspired from Ebin \& Marsden [EbM], the 
closure of the set $M_{c}$ in the $W^{s,2,\infty}$ topology: 
\begin{equation}
D^{s,2,\infty}_{c}(\Omega) \ = \ Cl^{W^{s,2,\infty}} \ M_{c} \ \ ,
\label{defgrup}
\end{equation} 
where $s$ is chosen to be greater than a critical value.

The tangent space at $D^{s,2,\infty}(\Omega)$ in $1_{\Omega}$ is denoted by 
$T_{1_{\Omega}} D^{s,2,\infty}(\Omega)$ and it is equal to
$$W^{s,2,\infty}_{0}(\Omega,R^{n}) \ = \ Cl^{W^{s,2,\infty}} \ \left\{ \eta \in
C^{\infty}(\Omega,R^{n}) \cap 
W^{s,2,\infty}(\Omega,R^{n}) \mbox{ : } supp \ \eta \subset \Omega \right\} \ \ .$$

 This tangent space may be endowed with the norm 
$\| \ \|_{1_{\Omega}}\ = \ \| \ \|_{W^{s,2,\infty}}$. The tangent space at $\phi$ and the
corresponding norm can be obtained from the tangent space and norm at the identity by right 
multiplication with $\phi$. 

The existence of the exponential map is proven in [EbM], theorem 3.1. for the group of 
$W^{s,2}$ diffeomorphisms of $M$, where $M$ is a compact manifold without boundary. The group
defined at (\ref{defgrup}) can be seen as a group of diffeomorphisms of a compact manifold equal to 
identity near a  point of the manifold. The results of Ebin \& Marsden are still true in this case. 
In theorem 3.4. from the
same article is proven that any continuous time dependent vector field can be integrated. 
We give the exact statements of these theorems, adapted to our
case. We do not give the proof, since it requires the same techniques as the ones employed in the  
article mentioned before: 

\vspace{.5cm}

{\bf Theorem 4.1.} {\it  Let $\Omega$ be an open bounded set in $R^{n}$ with piecewise 
smooth boundary and $s$ a number satisfying 
$$s \geq \frac{n}{2} + 2 \ \ \mbox{ .}$$ 

 Let $\eta \in W^{s,2,\infty}_{0}(\Omega,R^{n})$ be a vector field on $\Omega$, vanishing 
near $\paromega$. Then the one parameter flow generated by $\eta$ is a $C^{1}$ one parameter 
subgroup of $D^{s,2,\infty}_{c}(\Omega)$. 

Let $t \mapsto \eta_{t} \in  W^{s,2}_{0}(\Omega,R^{n})$ be a continuous time dependent vector 
field on  $\Omega$. Then the problem: 
$$ \eta_{t} = \dot{\phi}_{t} . \phi_{t}^{-1} \ \ \ , \phi_{0} = 1_{\Omega}$$ 
has a solution, unique, which is a $C^{1}$ curve in  $D^{s,2,\infty}_{c}(\Omega)$. }

\vspace{1.5cm}

We shall define further the notion of smooth crack propagation curve.

\vspace{.5cm}

{\bf Definition 4.1.} {\it Let us consider $s \geq \frac{n}{2} + 2$.  An 
$s$-smooth crack propagation curve  is a $C^{1}$ curve $$t \in [0,T] \mapsto \phi_{t} \in 
D^{s,2,\infty}_{c}(\Omega) \ \ ,$$ satisfying the following properties:
\begin{itemize}
\item
$\phi_{0} = 1_{\Omega}$, the identity map of $\Omega$, 
\item
$\phi_{t}(K) \subset \phi_{t'}(K)$ for any $0 \leq t < t' \leq T$. 
\end{itemize} }

\vspace{.5cm}

For any smooth crack propagation curve  we  associate the curve $t \mapsto K_{t} = \phi_{t}(K)$. The last curve has
the obvious property that for any $t > 0$ the crack set $K_{t}$ can be continuously deformed into $K$. In other 
words, $K_{t}$ and $K$ are topologically the same. We restrict our attention to this kind of evolution of the initial 
crack set. In this approach the initial crack $K$ may be as complex as we wish, but this complexity remains the same
during its evolution.

There are infinitely many smooth crack propagation curves $t \mapsto \phi_{t}$ 
with the same associated function $t \mapsto K_{t}$.

\vspace{.5cm}

{\bf Definition 4.2.} {\it Let $t \in [0,T] \mapsto \eu^{0}(t) \in C(\paromega,R^{n})$ be a 
$C^{1}$ curve of imposed
displacements on the exterior boundary of the body. An admissible fracture curve is any 
$C^{1}$ function 
$$t \in [0,T] \mapsto (\eu^{*}_{t},\phi_{t}) \in W^{1,2}_{K} \times D_{c}^{2,p,\infty}(\Omega) \
\ ,$$ 
satisfying the following items:

\begin{itemize}
\item
for any $t \in [0,T]$ $\eu^{*}_{t} = \eu^{0}(t)$ on $\paromega$,
\item
$t \mapsto \phi_{t}$ is a s-smooth crack propagation curve.
\end{itemize}  }

\vspace{.5cm}

Let us denote $\eu_{t} = \eu^{*}_{t} . \phi_{t}^{-1}$. Because for any $t$ $\phi_{t}$ equals the 
identity map near the boundary of $\Omega$, we deduce that $\eu_{t} = \eu^{0}(t)$ on $\paromega$. 
Remark also that the application $\eu^{*} \mapsto \eu^{*} . \phi^{-1}$ , where $\phi \in
D^{s,2,\infty}_{c}(\Omega)$ , maps $W^{1,2}_{K}$ onto $W^{1,2}_{\phi(K)}$, therefore $\eu_{t}$ is 
an admissible displacement of the body $\Omega$ with the crack set $\phi_{t}(K)$.

\vspace{.5cm}

{\bf Definition 4.3.} {\it 
A balanced fracture curve is an admissible fracture curve $t \mapsto (\eu^{*}_{t}, \phi_{t})$ such
that at any moment $t \geq 0$ the displacement $\eu_{t} = \eu^{*}_{t} . \phi_{t}^{-1}$ is a
solution of the problem }

$$min \ \left\{  \int_{\Omega} w(\nabla \eu) \mbox{ d}x \mbox{ : }  \eu \in 
W^{1,2,\infty}_{\phi_{t}(K)} (\eu^{0}(t)) \right\} \ \ .$$ 

\vspace{.5cm}

 For 
any smooth crack propagation curve there is an associated balanced fracture curve, unique to 
rigid displacements.

We are concerned now with the evolution of the area of the crack set. 

It is well-known (see [All]) that the variation of the area of the crack set $\phi_{t}(K)$ is 
$$\frac{d}{dt} \hen(\phi_{t}(K)) \ = \ \int_{\phi_{t}(K)} div_{s} \eta_{t} \mbox{ d}\hen \ \ ,$$ 
where $div_{s}$ is the tangential divergence with respect to the surface $\phi_{t}(K)$. If we
denote by $\en$ the field of normals to this surface, the expression of the tangential derivative 
$div_{s}$ is: 
$$ div_{s} \eta \ = \ \eta_{i,i} \ - \ \eta_{i,j} \en_{i} \en_{j} \ \ .$$

For a smooth crack propagation curve $t \mapsto \phi_{t}$ the condition that the crack 
grows implies that for any $t \geq 0$ 
\begin{equation}
\left\{ \begin{array}{ll}
\eta_{t} \cdot \en \ = 0 & \mbox{ on } \phi_{t}(K) \\
\int_{\phi_{t}(K)} div_{s} \eta_{t} \mbox{ d}\hen \ \geq 0 
\end{array} \right.
\label{difcond}
\end{equation}

It will be useful to define a perimeter measure of the edge of the crack. For this let us consider 
a crack set $K$. The perimeter of $\partial K$ can be defined, with the help of a flux-divergence
formula, like this: 
$${\cal P} (\partial K) \ = \ sup \left\{ \int_{K} div_{s} \eta \mbox{ d}\hen \mbox{ ; } 
\eta \in C^{\infty}_{0}(\Omega,R^{n}) \ , \mid \eta \mid \leq 1 \ ,  
\ \eta \cdot \en = 0 \mbox{ on } K \right\} \ \ .$$
The perimeter of $\partial K$ as a measure is defined first as an additive function over 
the set of open subsets $B$ of $\Omega$:
\begin{equation}
{\cal P} (\partial K) (B) \ = \ sup \left\{ \int_{K} div_{s} \eta \mbox{ d}\hen \mbox{ ; } 
\eta \in C^{\infty}_{0}(B,R^{n}) \ , \mid \eta \mid \leq 1 \ ,   
\ \eta \cdot \en = 0 \mbox{ on } K \right\}
\label{permes}
\end{equation}
This function generates a measure which is the perimeter of $\partial K$ as a measure.

\section{Crack propagation criterions}
\indent

In order to select one or more crack propagation curves among all admissible ones one needs 
a criterion of brittle crack propagation.   
There are several such  criterions. We shall discuss about two of them, 
in principle different. 

The Griffith criterion asserts that during the crack propagation the energy release rate due 
only to the crack evolution has to be greater than a critical quantity. This criterion may take 
different mathematical forms. Whatever this form may be, this criterion is formulated 
in terms of a critical speed. 

The Irwin type criterion asserts that during the crack propagation some intensity (or
concentration) factors have to be greater than a critical value. The original Irwin criterion 
is formulated in terms of stress intensity factors. However, it is straightforward that it can be 
reformulated by using  the elastic energy concentration on the edge of the crack. We shall refer to 
any crack propagation criterion using energy concentration factors (whatever they mean) as  
to the Irwin criterion.  

The goal of this section is to give precise mathematical expressions to the Griffith and Irwin
criterions for the class of smooth crack propagation curves.

\vspace{.5cm}

 We shall use in the sequel the notation 
$$ \eu(\phi_{t}, \eu^{0}(t)) \ = \ \eu_{t} \ = \  \eu(\phi_{t}(K), \eu^{0}(t))$$ 
since the initial crack $K$ is given. 
The stress field  associated to this displacement is $\sigma_{t} = \sigma (\nabla \eu_{t})$. 

The power communicated by the universe to the body at the moment $t$  has the expression 
$$P(t) \ = \ \int_{\paromega} \sigma_{t} \en \cdot \deu^{0}_{t} \mbox{ d}\hen \ \ .$$ 
The Griffith criterion of brittle fracture propagation is the following:

\vspace{.5cm} 

\hspace{2.cm} {\it A balanced fracture curve $t \mapsto  (\eu^{*}_{t}, \phi_{t})$ satisfies 
the Griffith criterion if at any moment $t \geq 0$ }

\begin{equation}
\frac{d}{dt} \left\{ \int_{\Omega} w(\nabla \eu(\phi_{t}, \eu^{0}(t))) \mbox{ d}x \ + \ 
G \hen(\phi_{t}(K)) \right\} \leq P(t) \ \ . 
\label{g1}
\end{equation}

\vspace{.5cm}
\noindent
$G$ is a material constant, named the constant of Griffith. 

The energy release rate due only to the propagation of the crack has the well-known expression:
\begin{equation}
{\cal E}(t) \ = \ P(t) \ - \  
\frac{d}{dt} \ \int_{\Omega} w(\nabla \eu(\phi_{t}, \eu^{0}(t))) \mbox{ d}x  \ \ .
\label{err}
\end{equation}

The inequality (\ref{g1}) is equivalent to 
$${\cal E}(t) \ \geq  \ G \ \frac{d}{dt}  \ \hen(\phi_{t}(K)) \ \ .$$

This criterion seems hard to work with it. The following proposition will lead us to an easier 
version of (\ref{g1}). 


\vspace{.5cm}

{\bf Proposition 5.1.} {\it Let $t \mapsto \phi_{t}$ be a smooth crack propagation curve, 
$\eta_{t} \ = \ \dot{ \phi_{t}} . \phi_{t}^{-1}$ and $\eu^{0} \in H^{\frac{1}{2}}(\paromega) 
\cap L^{\infty}(\paromega)$. Then for all $t \geq 0$ the following inequality is true : }
\begin{equation}
 \langle \frac{d}{dt} \left[  \eT(\phi_{t}) \right] \eu^{0} , \eu^{0} \rangle \ \leq  
\label{p1}
\end{equation}
$$ \leq  \int_{\Omega} \left\{ \left[ \eC \nabla \eu(\phi_{t},\eu^{0}) : \nabla \eu(\phi_{t},
\eu^{0}) \right] div \ \eta_{t} \ - 2 \left[ \eC \nabla \eu(\phi_{t},\eu^{0}) \right]_{ij} 
\left[ \nabla \eu(\phi_{t},\eu^{0}) \right]_{ik} \left[ \nabla \eta_{t} \right]_{kj} 
\right\} \mbox{ d}x \ \ .$$  

\vspace{.5cm}

{\bf Proof:} The smooth crack propagation curve $\phi_{t}$ is associated to the balanced fracture 
curve $( \eu(\phi_{t},\eu^{0}) . \phi_{t} , \phi_{t})$.  Let us choose a $t \geq 0$ and keep it 
fixed. We define now, for any $s \geq 0$,  
$$\eu_{t,s} = \eu(\phi_{t},\eu^{0}) . \phi_{t} . \phi_{t+s}^{-1} \ \ .$$
It is obvious that for all $s \geq 0$ $\eS_{\eu_{t,s}} \subset \phi_{t+s}(K)$, hence 
$$\forall s \geq 0 \ , \ \langle \eT(\phi_{t+s}) \eu^{0}, \eu^{0} \rangle \leq 
\int_{\Omega} \eC \nabla \eu_{t,s} : \nabla \eu_{t,s} \mbox{ d}x \ \ .$$
For $s = 0$ $\eu_{t,0} = \eu(\phi_{t},\eu^{0})$, so 
$$\langle \eT(\phi_{t}) \eu^{0}, \eu^{0} \rangle \ = \ 
\int_{\Omega} \eC \nabla \eu_{t,0} : \nabla \eu_{t,0} \mbox{ d}x \ \ .$$
We deduce that 
$$\langle \frac{d}{dt} \left[  \eT(\phi_{t}) \right] \eu^{0} , \eu^{0} \rangle = 
\frac{d}{ds} \langle \eT(\phi_{t+s}) \eu^{0}, \eu^{0} \rangle_{|_{s=0}} \leq 
\frac{d}{ds} \int_{\Omega} \eC \nabla \eu_{t,s} : \nabla \eu_{t,s} \mbox{ d}x_{|_{s=0}} \ \ .$$
It remains to prove that 
$$ \int_{\Omega} \left\{ \left[ \eC \nabla \eu(\phi_{t},\eu^{0}) : \nabla \eu(\phi_{t},
\eu^{0}) \right] div \ \eta_{t} \ - 2 \left[ \eC \nabla \eu(\phi_{t},\eu^{0}) \right]_{ij} 
\left[ \nabla \eu(\phi_{t},\eu^{0}) \right]_{ik} \left[ \nabla \eta_{t} \right]_{kj} 
\right\} \mbox{ d}x \ = $$ 
$$= \ \frac{d}{ds} \int_{\Omega} \eC \nabla \eu_{t,s} : \nabla \eu_{t,s} \mbox{ d}x_{|_{s=0}} \ \
.$$
The last equality results by direct calculation. It is sufficient  to make in the integral 
$$\int_{\Omega} \eC \nabla \eu_{t,s} : \nabla \eu_{t,s} \mbox{ d}x$$ the change of variables 
$y = \phi_{t} . \phi_{t+s}^{-1}(x)$ and then to perform the calculation of the derivative. 
\hfill $\bullet$

\vspace{.5cm} 

{\bf Definition 5.1.} {\it The generalized Rice's J integral is the functional }
$$K2 \ : \ D^{s,2,\infty}_{c}(\Omega) \times T_{1_{\Omega}}  D^{s,2,\infty}_{c}(\Omega) 
\times \left\{  H^{\frac{1}{2}}(\paromega) \ \cap L^{\infty}(\paromega) \right\} \ \rightarrow \ R 
\ \ ,$$
\begin{equation}
K2(\phi,\eta,\eu^{0}) \ = \ 
\int_{\Omega} \left\{ - \ \frac{1}{2} \left[ \eC \nabla \eu(\phi,\eu^{0}) : \nabla \eu(\phi,
\eu^{0}) \right] div \ \eta \ + \right.
\label{k2}
\end{equation}
$$\left.  +  \left[ \eC \nabla \eu(\phi,\eu^{0}) \right]_{ij} 
\left[ \nabla \eu(\phi,\eu^{0}) \right]_{ik} \left[ \nabla \eta \right]_{kj} 
\right\} \mbox{ d}x \ \ .$$  

\vspace{.5cm}

We shall explain why $K2$ is called "the generalized Rice's J integral" in the  remark 5.2.. Before
that   we give an alternative expression of $K2$. 

 Let us simplify the notations: 
$$\sigma = \eC \nabla \eu(\phi,\eu^{0}) \ \ , \  \eu = \eu(\phi,\eu^{0}) \ \ , \ w = \frac{1}{2} 
\eC \nabla \eu(\phi,\eu^{0}) : \nabla \eu(\phi,\eu^{0}) \ \ .$$ 
If $\eu$ is a $C^{2}$ function then the following equality makes sense: 
$$w \ \eta_{i,i} - \sigma_{ij} \eu_{i,k} \ \eta_{k,j} = 
\left[ w \ \eta_{i} \ - \ \sigma_{lj} \eu_{l,k} \
\eta_{k} \right]_{,i} - \sigma_{km} \eu_{k,mi} \ \eta_{i} + \sigma_{li} \eu_{l,ki}  \ \eta_{k} + 
\sigma_{li,i} \eu_{l,k}  \ \eta_{k} \ \ .$$ 
The divergence of the field $\sigma$ is equal to $0$, hence 
$$K2(\phi,\eta,\eu^{0}) \ = \ - \ \int_{\Omega} \left[w \ \eta_{i} \ - \ \sigma_{li} \eu_{l,k} \
\eta_{k} \right]_{,i} \mbox{ d}x \ \ .$$  
Let $\partial \phi(K)$ be the edge of the crack set $\phi(K)$. We define the tubular neighborhood
of $\partial \phi (K)$ of radius $r$:  
$$B_{r}= B_{r}(\partial \phi(K)) = \cup_{x \in \partial \phi(K)} B(x,r) \ \ .$$
The field of normals over $\partial \ B_{r}(\partial \phi(K))$ will be denoted by $\nu$, without 
specifying the parameter $r$. 

The function 
$$x \in \Omega \ \mapsto w \ div \eta (x) - \sigma_{li} \eu_{l,k} \ \eta_{k,i} \ (x)$$ 
is measurable, $\eta = 0$ on $\paromega$ and $\sigma \en = 0$ on $\phi(K)$, therefore 
$$K2(\phi, \eta,\eu^{0}) = \lim_{r \rightarrow 0} \left\{ \int_{\partial B_{r}(\partial \phi(K))} 
\left\{ - w \eta \cdot \nu \ + \sigma_{li} \eu_{l,k} \ \eta_{k} \ \nu_{i} \right\} \mbox{ d} \hen  
\right\} + \int_{\phi(K)}  \left[ w \right] \eta \cdot \en \mbox{ d}\hen \ \ .$$ 
If $\eta \cdot \en = 0$ on  $\phi(K)$ then 
\begin{equation}
 K2(\phi, \eta,\eu^{0}) = \lim_{r \rightarrow 0} \left\{ \int_{\partial B_{r}(\partial \phi(K))} 
 - \ \left\{ w \eta \cdot \nu \ - \sigma_{li} \eu_{l,k} \ \eta_{k} \ \nu_{i} \right\} \mbox{ d} \hen  
\right\} \ \ . 
\label{k2asrice}
\end{equation}

\vspace{.5cm}

{\bf Remark 5.1.}In [Oht1---4] Ohtsuka generalizes the Rice's integral J. For a given, 
very smooth, curve of increasing crack sets $t \mapsto K_{t}$, he founds a vector field 
$\eta_{0}$ whose one-parameter flow $\phi_{t}$ has the property 
$$K_{t} \ = \ \phi_{t}(K_{0}) \ \ .$$
He takes into consideration not only imposed displacements on the exterior boundary of the body but 
also imposed forces. The connection between his generalization and the functional $K2$ is 
the following: 
$$GJ(K_{t},\eu^{0}) \ = \ K2(\phi_{t},\eta_{0}, \eu^{0}) \ \ .$$
Ohtsuka  proves  that $GJ(K_{t},\eu^{0}(t))$ equals the energy release rate ${\cal E}(t)$. 
 \hfill $\bullet$

\vspace{.5cm}

{\bf Remark 5.2.}  In the
very particular case of one-dimensional crack propagation following a straight line Rice [R] 
expresses 
the energy release rate by an integral, named $J$. The expression of $J$ is formally similar with 
the right member of (\ref{k2asrice}) when $\eta$ is taken equal to a constant vector (the speed of 
the crack edge) parallel with the crack. 

The source of the notation $K2$ can be found in the 
particular case of a straight crack in the 
2-dimensional configuration. Let us suppose  that the body is a cylinder with section $\Omega$
and consider only anti-plane displacements: 
$$\eu: \Omega \times R \rightarrow R^{3} \ , \  \ \  \eu(x_{1},x_{2},x_{3}) \ = \ 
(0,0,u(x_{1},x_{2})) \ \ .$$ 
In this case the integral $J$ is proportional with the square of the stress intensity factor in
mode III denoted by $K_{III}$.   \hfill $\bullet$ 

\vspace{.5cm}

In the sequel we shall use the notation: 
$$\eT(\phi) \ = \ \eT(\phi(K)) \ \ .$$ 
The energy release rate ${\cal E}(t)$ has the following expression: 
$${\cal E}(t) \ =  \ \langle \eT\phi_{t}) \eu^{0}(t) , \deu^{0}(t) \rangle \ - \ 
\frac{1}{2} \ \frac{d}{dt} \left\{ \langle  \eT\phi_{t}) \eu^{0}(t) , \eu^{0}(t) \rangle \right\} 
\ \ .$$ 
We use the fact that $\eT(\phi)$ is self-adjoint for proving that 
$$\langle \eT\phi_{t}) \eu^{0}(t) , \deu^{0}(t) \rangle \ = \ 
\langle \eT\phi_{t}) \deu^{0}(t) , \eu^{0}(t) \rangle  \ \ .$$
The last two relations show  that the energy release rate has the expression: 
$$ {\cal E}(t) \  = \ - \  
\frac{1}{2} \langle \ \frac{d}{dt} \left[  \eT(\phi_{t}) \right] \eu^{0}(t) , \eu^{0}(t) \rangle
\  \ . $$ 
We  use  proposition 5.1. and definition  5.1. for proving that
\begin{equation} 
{\cal E}(t) \  \geq \   K2(\phi_{t}, \eta_{t}, \eu^{0}(t)) \ \ .
\label{errk2}
\end{equation}

We propose now a Griffith type criterion of smooth fracture propagation.

\vspace{.5cm}

\hspace{2.cm} {\it A smooth crack propagation curve $t \mapsto \phi_{t}$ satisfies the generalized 
Griffith criterion if at any moment $t \geq 0$} 
\begin{equation}
K2(\phi_{t},\eta_{t},\eu^{0}(t)) \ \geq \ G \ 
\| \eta_{t} \|_{\infty} {\cal P}(\partial \phi_{t}(K)) (supp \ \eta_{t})
\label{g2}
\end{equation}

\vspace{.5cm}

In our criterion of brittle fracture propagation ${\cal E}(t)$ is replaced by 
$K2(\phi_{t}, \eta_{t}, \eu^{0}(t))$ and the variation of the area of the crack 
$$\int_{\phi_{t}(K)} div_{s} \eta_{t} \mbox{ d}\hen$$ is replaced by 
$\| \eta_{t} \|_{\infty} {\cal P}(\partial \phi_{t}(K)) (supp \ \eta)$. 
This criterion is therefore stronger 
 than the classical Griffith criterion (this is a consequence of (\ref{errk2}) and of the
definition of the perimeter (\ref{permes})). 
The  two selection criterions are however equivalent if at any 
moment $t$ the crack is smooth enough and it evolves in a very smooth manner.  

\vspace{1.5cm}

The Irwin type criterion of brittle crack propagation is formulated in terms of concentration 
coefficients of the elastic energy. 

\vspace{.5cm}

{\bf Definition 5.2.} {\it Let $\eu$ be a special function with bounded variation.  The lower 
(respectively upper) concentration coefficients of the elastic energy  of $\eu$ are }
\begin{equation}
C2^{-}(\eu,x) = \frac{1}{\pi} \ \liminf_{r \rightarrow 0} 
\frac{\int_{\Omega \cap B(x,r)} w(\nabla \eu) \mbox{ d}x}{r^{n-1}} \ \ ,
\label{lower} 
\end{equation}
\begin{equation}
C2^{+}(\eu,x) = \frac{1}{\pi} \  \limsup_{r \rightarrow 0} 
\frac{\int_{\Omega \cap B(x,r)} w(\nabla \eu) \mbox{ d}x}{r^{n-1}} \ \ .
\label{upper}
\end{equation}

\vspace{.5cm}

If $x \in \Omega \setminus \eS_{\eu}$ then $C2^{+}(\eu,x) = 0$, because $w(\nabla \eu)$ has 
approximate limit in $x$. 

Let us denote by $E\eu$ the symmetric part of the distributional derivative of $\eu$: 
$$E \eu \ = \ \frac{1}{2} \left( D \eu \ + \ D^{T}\eu \right) \ \ .$$
$E\eu$ can be seen, like $D\eu$, as a vector measure; let $\mid E\eu \mid$ be the scalar measure 
of the variation of $E\eu$. For any borelian set $B \subset \Omega$ 
\begin{displaymath}
\mid E\eu \mid  (B) \ = \  sup \ \left\{ \sum^{\infty}_{i=1} \mid E\eu  (A_{i}) \mid 
\mbox{ : } \cup_{i=1}^{\infty} A_{i} \subset B \ , \ A_{i} \cap A_{j} = \emptyset \ \ \forall 
i \not = j \right\} \ \ .
\end{displaymath}
Let us consider the set $\Theta_{\eu}$ introduced by Kohn [K]: 
$$\Theta_{\eu} = \left\{ x \in \Omega \mbox{ : } \limsup_{r \rightarrow 0} 
\frac{ \mid E\eu \mid (B(x,r))}{r^{n-1}} \ > \ 0 \right\}$$
A result from [ACDM] ,concerning special functions with bounded deformations in the
particular case of $\eu$ special function with bounded variation, assures us that $\Theta_{\eu}$
differs from $\eS_{\eu}$ by a set which has $\hen$ measure zero. We use the growth condition of 
$w$ (\ref{grow}) in order to prove that if $x \in \Theta_{\eu}$ then $C2^{+}(\eu,x) = 0$. 
Therefore we have the following localization of the support of the energy concentration 
coefficients:
$$supp \ C2^{-}(\eu, \cdot) \ \subset \ supp \ C2^{+}(\eu, \cdot) \subset \eS_{\eu} \setminus 
\Theta_{\eu} \ \ .$$
One can prove, using the same proof, that if $K$ is a crack set and $\eu \in W^{1,2}_{K} \cap 
L^{\infty}(\Omega,R^{n})$ then 
$$supp \ C2^{-}(\eu, \cdot) \ \subset \ supp \ C2^{+}(\eu, \cdot) \subset \partial K \ \ .$$

The Irwin type criterion that we propose is the following: 

\vspace{.5cm}

\hspace{2.cm} {\it A smooth crack propagation curve $t \mapsto \phi_{t}$ satisfies the 
Irwin type criterion if at any moment $t \geq 0$ and for any $x \in  \ \partial  \phi_{t}(K)$ } 
\begin{equation}
C2^{+}(\eu(\phi_{t},\eu^{0}(t)) , x) \ \geq \ G \ \ . 
\label{irwin}
\end{equation}

\vspace{.5cm}

\section{Energy concentration and crack propagation}
\indent

The goal of this section is to find the connections between the Griffith criterion (\ref{g2}) and 
the Irwin criterion (\ref{irwin}), using a measure-theoretic approach. 

The first step is to associate a measure to the generalized Rice integral $K2$. 

\vspace{.5cm}

{\bf Definition 6.1.} {\it Let us consider the linear functional 
$$k2 \ : \ \eSBV(\Omega,R^{n}) \cap 
L^{\infty}(\Omega,R^{n}) \ \times \ C^{\infty}_{0}(\Omega,R^{n}) \ \rightarrow R \ \ ,$$
\begin{equation}
k2(\eu,\eta) \ = \ \int_{\Omega} \left\{ \sigma(\nabla \eu)_{ij} \eu_{i,k} \eta_{k,j} \ - \ w(\nabla \eu) 
\ div \ \eta \right\} \mbox{ d}x \ \ .
\label{k2pre}
\end{equation}
Then, for any $\eu \in \eSBV(\Omega,R^{n}) \cap 
L^{\infty}(\Omega,R^{n})$ and for any open subset $B$ of $\Omega$ we define 
\begin{equation}
\mid K2 \mid (\eu) (B) \ = \ sup \left\{ k2(\eu,\eta) \mbox{ : } supp \ \eta \subset \subset B \ , \ 
\mid \eta \mid \leq 1 \ , \  \eta \cdot \en = 0 \mbox{ on } \eS_{\eu} \right\} \ \ .
\label{k2mes}
\end{equation}  
Since $\mid K2 \mid(\eu)( \cdot)$ is finitely additive, it gives raise to a positive measure,
denoted by $\mid K2 \mid(\eu)$ too, named the generalized Rice integral as a measure. }

\vspace{.5cm}

It is easy to see that for a smooth crack propagation curve $t \mapsto \phi_{t}$ it is true that 
$$\forall \ t \geq 0 \ , \ K2(\phi_{t}, \eta_{t},\eu^{0}(t)) \ = \ k2(\eu(\phi_{t},\eu^{0}(t)) , 
\eta_{t}) \ \leq \ \mid K2 \mid (\eu(\phi_{t},\eu^{0}(t))) (supp \ \eta_{t} )$$

\vspace{.5cm}

{\bf Remark 6.1.} The measures ${\cal P}(K)$ and $\mid K2 \mid(\eu)$ are constructed in the 
same way. The departure point was in both cases a linear functional over $C^{\infty}_{0}(\Omega, 
R^{n})$ --- the area variation integral in the first case and the $k2$ integral in the second case. 
\hfill $\bullet$

\vspace{.5cm}

We shall define now the upper energy concentration coefficient as a measure. 

\vspace{.5cm}

{\bf Definition 6.2.} {\it Let us consider $\eu \in \eSBV(\Omega,R^{n})$. For any  $r > 0$ 
let  $B_{r}(\eu)$ be a tubular neighborhood of $\eS_{\eu}$ of radius $r$. For any 
open set $B$ in $\Omega$ we define the upper energy concentration coefficient: 
\begin{equation}
CM^{+}(\eu)( B) = \limsup_{r \rightarrow 0} \frac{ \int_{B_{r}(\eu) \cap B} w(\nabla \eu) 
\mbox{ d}x}{r} \ \ .
\label{cm+}
\end{equation}
Since $CM^{+}(\eu)(\cdot)$ is an additive function, it give raise to a positive measure denoted by 
$CM^{+}(\eu)$ too, named the upper energy concentration coefficient as a measure. }

\vspace{.5cm}

$CM^{+}(\eu)$  splits in two parts: the absolute continuous part with respect to $\leb$ and the 
singular one. 
$$CM^{+}(\eu) \ = \ CM^{+}(\eu)^{A} \ + \ CM^{+}(\eu)^{S} \ \ .$$
Let us consider, for $\eu \in \eSBV(\Omega,R^{n})$, the decomposition of $D\eu$: 
$$D\eu \ = \ \nabla \eu \mbox{ d}x \ + \ [\eu] \otimes \en \mbox{ d}\hen_{|_{\eS_{\eu}}} \ \ .$$
The measure $ CM^{+}(\eu)^{S}$ can be decomposed in two parts --- the absolute continuous part 
with respect to  
the measure $$J(\eu)\ = \  [\eu] \otimes \en \mbox{ d}\hen_{|_{\eS_{\eu}}}$$
and the remaining singular part --- i.e. :
$$CM^{+}(\eu)^{S} \ = \  CM^{+}(\eu)^{J} \ + \ CM^{+}(\eu)^{C} \ \ .$$
It is easy to see that  $CM^{+}(\eu)^{A} \ = \ 0$. 
We obtain therefore the following decomposition of  $CM^{+}(\eu)$:
\begin{equation}
CM^{+}(\eu) \ =  \ CM^{+}(\eu)^{J} \ + \ CM^{+}(\eu)^{C}
\label{firstdec}
\end{equation}
 If $\eu \in   W^{2,2}_{loc}(\Omega \setminus
\eS_{\eu},R^{n})$ and $\eS_{\eu}$ is a crack set, then we find, by direct calculation, the 
expression of $CM^{+}(\eu)^{C}$: 
\begin{equation}
 CM^{+}(\eu)^{C} \ = \  C2^{+}(\eu,x) \mbox{ d}{\cal H}^{n-2}_{|_{\partial \eS_{\eu}}} \ \ . 
\label{cantorconc}
\end{equation}

\vspace{.5cm}

The main result of the paper is the following: 

\vspace{.5cm}

{\bf Theorem 6.1.} {\it  Let $K$ be a crack set and $\eu \in W^{1,2,\infty}_{K}(\eu^{0}) \cap
W^{2,2}_{loc}(\Omega \setminus K,R^{n})$. Then $\eu = \eu(K,\eu^{0})$ if and only if  for 
any open set $B \subset \Omega$ }
\begin{equation}
\mid K2 \mid (\eu) (B) \ \leq \ CM^{+}(\eu)^{C}(B) \ \ . 
\label{girv}
\end{equation}

\vspace{.5cm}

A consequence of theorem 6.1. and (\ref{girv}) is that the Griffith criterion (\ref{g2})  is 
stronger than the Irwin type criterion (\ref{irwin}). Indeed, let us consider a smooth crack
propagation curve $t \mapsto \phi_{t}$ satisfying the generalized Griffith criterion 
(\ref{g2}). Then, at any moment $t \geq 0$, we have: 
$$\mid K2 \mid (\eu(\phi_{t}, \eu^{0}(t)))(\Omega) \  \| \eta_{t} \|_{\infty} \ \geq \ K2(\phi_{t},\eta_{t},\eu^{0}(t)) \ 
\geq \ \ G \ 
\| \eta_{t} \|_{\infty} {\cal P}(\partial \phi_{t}(K)) (supp \ \eta_{t}) \ \ .$$
This chain of inequalities, (\ref{girv}) and (\ref{cantorconc}) imply the following relation: 
$$\int_{\partial \phi_{t}(K)} C2^{+}(\eu (\phi_{t}, \eu^{0}(t)), \cdot) \mbox{ d}\hen \ \geq  \
G \ {\cal P}(\partial \phi_{t}(K)) (supp \ \eta_{t}) \ \ .$$
Therefore at any moment $t \geq 0$ there is at least an $x \in \partial \phi_{t}(K)$ such that 
$$C2^{+}(\eu (\phi_{t}, \eu^{0}(t)), x) \ \geq G \ \ .$$

\vspace{.5cm}

\section{Proof of theorem 6.1.}
\indent

 {\it $1^{st}$ implication.}
 We want to prove that if $\eu = \eu(K,\eu^{0})$ then for any 
open set $B$ the relation (\ref{girv}) is true, namely 
$$\mid K2 \mid (\eu) (B) \ \leq \ CM^{+}(\eu)^{C}(B) \ \ .$$
 
 Let us consider a vector field  $\eta$, such that: 
\begin{equation}
\begin{array}{ll} 
\eta \in C^{\infty}_{0}(\Omega,R^{n}) &   \\
\mid \eta \mid \leq 1 &  \\
\eta \cdot \en = 0 & \mbox{ on } K  \\
 supp \ \eta \subset B & .
\end{array}  
\label{preta}
\end{equation}
Any such vector field $\eta$ generates an one-parameter flow $s  \mapsto \phi_{s}$. 
We can always find a curvilinear 
coordinate system $(\alpha_{1}, . . . , \alpha_{n-1}, \gamma)$ such that on the edge of the crack
set $\partial K$ $\gamma = 0$ ,  the surface $\gamma = s$ is the boundary  of a tubular 
neighborhood of $\partial K$, namely  
$$B_{r(s)} (\partial K) \ = \ \cup_{x \in \partial K} B_{r(s)} (x)$$
and for all $s > 0$ 
\begin{equation}
\phi_{s}(K) \subset \ K \cup B_{r(s)} (\partial K) \ \ .
\label{subuni}
\end{equation} 
We can suppose moreover that 
$$\lim_{s \rightarrow 0} \frac{r(s)}{s} \ = \ 1$$
because (\ref{preta}) affirms that the crack $\phi_{s}(K)$ evolves with sub-unitary speed. 

We denote by $\eu_{s}$  the solution of the problem:
\begin{equation}
\left\{ \begin{array}{ll}
div \ \sigma(\nabla \eu) \ = \ 0 & \mbox{ in } \Omega \setminus \left( K
\cup  B_{r(s)} \right) \\
 \sigma(\nabla \eu) \en \ = \ 0 & \mbox{ on } \partial B_{r(s)} (\partial K) \\
\sigma(\nabla \eu) \en \ = \ 0 & \mbox{ on (both sides of) } K \setminus B_{r(s)} (\partial K) \\
\eu \ = \ \eu^{0} & \mbox{ on } \partial \Omega  \ \ . 
\end{array}
\right.
\label{eus}
\end{equation}
We consider, for any $s > 0$, the following  stress field: 
\begin{equation}
\sigma_{s}(x) \ = \ \left\{ \begin{array}{ll}
  \sigma(\nabla \eu_{s})(x) & \mbox{ if } x \in  \Omega \setminus \left( K
\cup  B_{r(s)} \right) \\
0 & \mbox{ elsewhere } \ \ .  
\end{array} 
\right. 
\end{equation}
It is obvious that for any $s$  $\sigma_{s}$ is admissible with respect to $\phi_{s}(K)$ 
(definition 3.1.1).  We use (\ref{moreau}) and proposition 5.1. to  deduce that 

\begin{equation} 
\liminf_{s \rightarrow 0} \frac{ \int_{\Omega} w(\nabla \eu_{s}) \mbox{ d}x \ - \ \int_{\Omega} w(\nabla \eu) \mbox{ d}x}{s} \ 
\leq \ 
\label{apulim}
\end{equation}
$$ \leq \ \lim_{s \rightarrow 0} 
\frac{ \int_{\Omega} w(\nabla \eu(\phi_{s},\eu^{0})) \mbox{ d}x \ - \ 
\int_{\Omega} w(\nabla \eu) \mbox{ d}x}{s} \leq \ - \ K2(1_{\Omega}, \eta, \eu^{0}) \ \ .$$
Our goal is  to prove that 
\begin{equation}
\liminf_{s \rightarrow 0} \frac{ \int_{\Omega} w(\nabla \eu_{s}) \mbox{ d}x \ - \ \int_{\Omega} 
w(\nabla \eu) \mbox{ d}x}{s} \ = \ 
 - \ \limsup_{s \rightarrow 0} \frac{ \int_{B_{r(s)} \cap B} w(\nabla \eu) \mbox{ d}x}{r(s)}  \ \ . 
\label{pulim}
\end{equation}
For this it is sufficient to show that 
\begin{equation}
\lim_{s \rightarrow 0} \frac{ \int_{\Omega \setminus B_{r(s)} } w(\nabla \eu_{s}) \mbox{ d}x \ - \ 
\int_{\Omega \setminus B_{r(s)}} 
w(\nabla \eu) \mbox{ d}x}{s} \ = \ 0 \ \ . 
\label{nec}
\end{equation}
Indeed, if (\ref{nec}) is true then an easy calculation shows that (\ref{pulim}) is true. With 
this in mind we return to (\ref{apulim}) and we see that 
$$- \ \limsup_{s \rightarrow 0} \frac{ \int_{B_{r(s)} \cap B} w(\nabla \eu) \mbox{ d}x}{r(s)}  \
\leq \ -K2(1_{\Omega},\eta,\eu^{0}) \ \ ,$$ 
hence 
$$K2(1_{\Omega},\eta,\eu^{0}) \ \leq \ \limsup_{s \rightarrow 0} 
\frac{ \int_{B_{r(s)} \cap B} w(\nabla \eu) \mbox{ d}x}{r(s)} \ \ .$$ 
We  take the supremum of the left side term  of the last inequality, over all 
$\eta$ satisfying  (\ref{preta}), and we obtain the inequality 
$$\mid K2 \mid (\eu) (B) \ \leq \ 
\limsup_{s \rightarrow 0} \frac{ \int_{B_{r(s)} \cap B} w(\nabla \eu) \mbox{ d}x}{r(s)} \ \ .$$ 
This inequality is equivalent with (\ref{girv}) (as it is shown by definition 5.2. and 
(\ref{cantorconc}) ). 

\vspace{.5cm}

{\it Proof of } (\ref{nec}). We consider for any $s$ the function 
$$\tilde{\eu_{s}}(x)  \ = \ \left\{ \begin{array}{ll} 
\eu_{s}(x) & \mbox{ if } x \not \in  B_{r(s)}(\partial K) \\
0 & \mbox{ if } x \in B_{r(s)}(\partial K)  \ \ \ . 
\end{array} \right. $$  
Each $\tilde{\eu_{s}}$ is a special function with bounded variation 
on $\Omega$.

In the following we shall proceed as in [BBC]. We denote by $\eve_{1}, ... ,
\eve_{n-1}, \eve_{\gamma}$ the local basis associated to the previous mentioned system of
coordinates. $\eve^{1}, ... ,
\eve^{n-1}, \eve^{\gamma}$ represents the dual local basis. W can choose the system of coordinates 
such that $\eve_{\gamma} = \ h \ \en$, where $\en$ is the normal of the surface $\gamma = ct.$ and 
$h$ a scalar function. 

The expression of the divergence of a tensor field in this system of coordinates is 
$$div \ \sigma \ = \ \left( \nabla_{\eve_{\alpha_{i}}} \ \sigma \right) \cdot \eve_{\alpha_{i}} \ +  \ 
\left( \nabla_{\eve_{\gamma}} \ \sigma \right) \cdot \eve_{\gamma} \ \ .$$  
We make the notation 
$$div_{\tau} \ \sigma \ = \ \left( \nabla_{\eve_{\alpha_{i}}} \ \sigma \right) \cdot \eve_{\alpha_{i}}$$
and we remark that in $div_{\tau} \ \sigma$ enters only partial derivatives with respect to 
$\alpha_{i}$. The divergence of the field $\sigma$ can be written as
$$div \ \sigma \ = \ div_{\tau} \ \sigma \ + \frac{1}{h} \frac{\partial}{\partial \gamma} \left( \sigma 
\cdot \en \right) \ - \ \frac{1}{h} \sigma \cdot \frac{\partial \en}{\partial \gamma} \ \ .$$

Consider $\sigma = \sigma (\nabla \eu)$ and 
$$\sigma^{m} (\alpha_{i}, s) \ = \ \sigma_{s}(\alpha_{i}, s) \ \ .$$ 
The field $\sigma^{m}$ has the following interpretation: imagine that at the "moment" $s$ the 
curvilinear cylinder $B_{r(s)}(\partial K)$ is removed from the elastic body $\Omega$. Then 
$\sigma^{m} (\cdot, s)$ represents the superficial stress of the remaining body on the new-created 
surface, i.e.  $\partial B_{r(s)}(\partial K)$. Therefore 
 $$\sigma^{m}(\alpha_{i},a) \en \ = \ 0 \ \ ,$$ 
where $\en$ is the normal of the surface $\gamma = a$.  

We consider, for any $s$ the auxiliary elastic problem: 
 \begin{equation}
\left\{ \begin{array}{lr}
 div \ \sigma^{aux}_{s} \ = \ 0 & \mbox{ in }  \Omega \setminus \left( K \cup B_{r(s)} (\partial K) \right) \\
 \sigma^{aux}_{s} \en  \ = \ \sigma(\alpha_{i}, s) \en & \mbox{ on } \left( K \setminus B_{r(s)}(\partial K) \right) \cup 
\partial B_{r(s)}(\partial K) \\
 \eu^{aux}_{s} \ = \ 0 & \mbox{ on } \paromega
\end{array} \right.
\label{paux}
\end{equation}
with the solution $\eu^{aux}_{s}(\alpha_{i}, \gamma)$, 
$\gamma \leq s$. The associated stress field is $\sigma^{aux}_{s}(\alpha_{i},\gamma)$.  
The superposition principle applied to the elastic problem (\ref{diric}) with the solution $\eu = \eu(K,\eu^{0})$, 
(\ref{eus}) and (\ref{paux}) leads us to the equality 
\begin{equation}
\sigma(\alpha_{i}, s ) \ = \ \sigma^{m}(\alpha_{i}, s) + \sigma^{aux}_{s}(\alpha_{i}, s) \ \ .
\label{deco}
\end{equation}
The following equality is obtained from a linear combination of the divergence equations 
contained in the problems (\ref{diric}), 
(\ref{eus}) and (\ref{paux}) , written in curvilinear
coordinates, together with (\ref{deco}).   
\begin{equation}
\frac{\partial}{\partial  \gamma} \left( \sigma \en \right) (\alpha_{i}, \gamma = s ) \ = \ 
 \frac{\partial}{\partial  \gamma}\left( \sigma^{aux}_{s} \en \right) (\alpha_{i},\gamma = s ) \ - \
h \ div \ \sigma^{m} (\alpha_{i},\gamma = s ) \ \ .
\label{inter}
\end{equation}
We make the derivative with respect to $s$ in the equality (\ref{deco}), then we multiply by $\en$ and we obtain 
\begin{equation}
\frac{\partial}{\partial  \gamma} \left( \sigma \en \right) (\alpha_{i}, \gamma = s ) \ = \ 
\frac{\partial}{\partial  \gamma'}_{|_{\gamma' = s}} \left( \sigma^{aux}_{\gamma'} \en \right) 
(\alpha_{i}, \gamma = s ) \ + \  \frac{\partial}{\partial  \gamma} 
\left( \sigma^{aux}_{s} \en \right) (\alpha_{i},\gamma = s) \ .  
\label{inter1}
\end{equation}
(\ref{inter}) and (\ref{inter1}) lead us to the equality 
\begin{equation}
\frac{\partial}{\partial  \gamma'}_{|_{\gamma' = s}} \left( \sigma^{aux}_{\gamma'} \en \right) 
(\alpha_{i}, \gamma = s ) \ = \ - \ h \ div \ \sigma^{m} (\alpha_{i},\gamma = s ) \ \ .
\label{clim}
\end{equation}

We keep now $s$ constant. Let us consider $s' \in [0,s]$. We use the notation
$$\Omega_{s'} = \Omega \setminus \left( K \cup B_{r(s')}(\partial K) \right) \ \ ,$$
$$S_{s'} = \partial B_{r(s')}(\partial K) \ \ .$$
It is obvious that if $s' < s"$ then $\Omega_{s"} \subset \Omega_{s'}$.

The field  $\sigma^{aux}_{s'}$ is by definition the elastic stress field in the body $\Omega_{s'}$, 
with the imposed displacement on it's exterior boundary $\paromega$ equal to $0$ and   
subjected to the loads $\sigma(\alpha_{i}, r( s')) \en(\alpha_{i}, s')$ on it's interior boundary 
$S_{s'}$. It is straightforward that the same field represents the elastic stress in the body 
$\Omega_{s}$, with the imposed displacement  on it's exterior boundary $\paromega$ equal to $0$ and 
 subjected to the loads $\sigma^{aux}_{s'}(\alpha_{i}, s) \en(\alpha_{i}, s)$ on $S_{s}$. 
Therefore, by derivation with respect to $s'$, we can prove that $\partial / \partial 
\gamma'_{|_{\gamma' = s}} \left( \sigma^{aux}_{s'} \right)$ is the elastic stress field 
in $\Omega_{s}$, resulting from the imposed displacement $0$ on $\paromega$  and the loads 
$\partial / \partial 
\gamma'_{|_{\gamma' = s}} \left( \sigma^{aux}_{s'} \en \right)(\alpha_{i},s)$ applied on $S_{s}$. 
The equation (\ref{clim}) shows that these loads are equal to 
$- \ h \ div \ \sigma^{m} (\alpha_{i}, s )$.

By our choice of the coordinate system we have the equality $h \ = \ \partial \ r(s)/ \partial s$. 
For our purposes it is not restrictive to assume that when $s$ tends to $0$ $\eu_{s}$ tends to 
$\eu$ and $\varepsilon ( \eu_{s})$ tends to $\varepsilon ( \eu)$. Indeed, 
there is  a positive constant $M$ such that for all $s \leq 1$  
$$ c \ \int_{\Omega} \mid \nabla \tilde{\eu_{s}} \mid^{2} \mbox{ d}x + \hen(\eS_{\tilde{\eu_{s}}}) 
\ \leq 
\  \int_{\Omega} w(\varepsilon( \tilde{\eu_{s}})) \mbox{ d}x + \hen(\eS_{\tilde{\eu_{s}}})  \ \leq \ M$$
where $c$ is the ellipticity constant of the tensor $\eC$. The last 
inequality is true because for all 
$s$ the elastic energy of $\tilde{\eu_{s}}$ is smaller that the elastic energy of $\eu$ and 
$\eS_{\tilde{\eu_{s}}} \subset \ K \cup \partial B_{r(s)}(\partial K)$. 
By the compactness theorem of Bellettini, Coscia \& Dal Maso [BCDM] 
and the regularity of $\tilde{\eu_{s}}$,  
there exists a sequence $s_{h} \rightarrow 0$ such that 
$\tilde{\eu_{s_{h}}}$ converges punctually to $\eu$ and $\varepsilon( \tilde{\eu_{s_{h}}})$ 
converges punctually to 
$\varepsilon( \eu)$.  

 Therefore when $s$ tends to $0$ 
$\sigma^{m}(\alpha_{i},s) = \sigma_{s}(\alpha_{i},s)$ tends to $\sigma(\alpha_{i},s)$, 
hence  $div \ \sigma^{m}$ tends to $0$ (because $div  \ \sigma = 0$). The conclusion is that  
$- \ h \ div \ \sigma^{m} (\alpha_{i}, s )$ tends to $0$ too, so the displacement solution of the 
problem 
$$\left\{ \begin{array}{ll}
div \ \sigma(\nabla \ev) \ = \ 0 & \mbox{ in }  \Omega_{s} \\
\ev \ = \ 0 & \mbox{ on } \paromega \\ 
\sigma (\nabla \ev) \en \ = \ 0 & \mbox{ on (both sides of) } K \cap \Omega_{s} \\ 
\sigma (\nabla \ev) \en \ = \ - \ h \ div \ \sigma^{m} (\alpha_{i}, s ) & \mbox{ on } \eS_{s} \ \ , 
\end{array} \right. $$
which is $\partial / \partial 
\gamma'_{|_{\gamma' = s}} \left( \eu^{aux}_{s'} \right)$, tends to $0$. In particular 
\begin{equation} 
\lim_{s \rightarrow 0} \frac{ \eu_{s} \ - \ \eu}{s} \ = \ - \ 
\frac{\partial}{\partial \gamma'}_{|_{\gamma' = 0}} \eu^{aux}_{\gamma'} \ = \ 0 \ \ .
\label{ulim}
\end{equation}

We remark that  
$$\lim_{s \rightarrow 0} \frac{ \int_{\Omega \setminus B_{r(s)} } w(\nabla \eu_{s}) \mbox{ d}x \ - \ 
\int_{\Omega \setminus B_{r(s)}} 
w(\nabla \eu) \mbox{ d}x}{s} \ = \ \frac{1}{2} \lim_{ s \rightarrow 0} 
\int_{S_{s}} \left( \sigma \en \right) \cdot  \frac{ \eu \ - \ \eu_{s}}{s} \mbox{ x} \hen \ \ .$$
This equality together with (\ref{ulim}) concludes the proof of (\ref{nec}) \hfill $\bullet$

\vspace{2.cm}

{\it $2^{nd}$ implication.}  We want to prove now that if $$\eu \in W^{1,2,\infty}_{K}(\eu^{0}) \cap
W^{2,2}_{loc}(\Omega \setminus K,R^{n})$$ and   for 
any open set $B \subset \Omega$ 
$$\mid K2 \mid (\eu) (B) \ \leq \ CM^{+}(\eu)^{C}(B) \ \ ,$$ 
then $\eu \ = \ \eu(K,\eu^{0})$. 

The regularity assumptions concerning $\eu$ and $K$ allow us to make the further calculations. 
As previously, we use the notations: $\sigma = \sigma(\nabla \eu)$, $w \ = \ w(\nabla \eu)$. 

Let us consider a vector field $\eta$ satisfying (\ref{preta}). Then 
$$k2(\eu,\eta) \ = \ \sigma_{ij} \eu_{i,k} \eta_{k,j} \ - \ w \ div \ \eta \mbox{ d}x \ = \ 
\int_{\Omega} - \sigma_{li,i} \eu_{l,k} \ \eta_{k} \ - \ \left[ w \ \eta_{i} \ - \ \sigma_{lj} 
\eu_{l,k} \eta_{k} \right]_{,i} \mbox{ d}x \ \ .$$
We deduce that the absolute continuous part of the measure $\mid K2 \mid (\eu)$ has the density 
$\mid \sigma_{li,i} \eu_{l,k} \mid$. The hypothesis (\ref{girv}) and the decomposition 
(\ref{firstdec}) 
 imply that this density is 
equal to $0$, hence: 
$$ div \ \sigma(\nabla \eu) \ = \ 0 \ \ \ \ \mbox{a.e. in } \Omega \ \ .$$
Therefore $\sigma$ is an admissible stress in the sense of the definition 3.1.1. . The before 
mentioned result of Del Piero [DP] implies that $\sigma \en$ has no jumps over $K$. 

The measure $\mid K2 \mid$ has a part concentrated on $K$, i.e. the absolute continuous 
part of  $\mid K2 \mid$ with respect to the measure $\hen_{|_{K}}$. The density of this part  
can be calculated from the expression of $k2$: 
$$k2(\eu,\eta) \ = \ \int_{\Omega} - \ \left[ w \ \eta_{i} \ - \ \sigma_{lj} 
\eu_{l,k} \eta_{k} \right]_{,i} \mbox{ d}x$$
Indeed, (\ref{preta}) and the regularity assumptions over $K$ and $\eu$ allow us to use a
flux-divergence formula in order to prove that the 
density of this part of $\mid K2 \mid$ is $\mid \sigma \en \mid$. Again, 
(\ref{girv}) and the decomposition 
(\ref{firstdec}) 
 imply that this density is 
equal to $0$, hence: 
$$\sigma \en \ = \ 0 \ \ \ \ \mbox{ on (both sides of) } K$$
This concludes the proof. \hfill $\bullet$

\section{Final considerations}
\indent

Let us return to relation (\ref{apulim}) from the proof of theorem 6.1., first implication.
 We have 
considered there $\eu = \eu(K,\eu^{0})$ and $\eta$ an arbitrary vector field, such that 
\begin{equation}
\begin{array}{ll} 
\eta \in C^{\infty}_{0}(\Omega,R^{n}) &   \\
\mid \eta \mid \leq 1 &  \\
\eta \cdot \en = 0 & \mbox{ on } K \ \ . 
\end{array}  
\label{didi}
\end{equation} 
This field generates the one-parameter flow $\phi_{s}$, which represents an arbitrary smooth crack propagation curve with sub-unitary initial velocity 
$\eta$. 
The energy release rate, defined by (\ref{err}), associated to this curve, at the moment $s=0$ will be denoted by ${\cal E}(\eta)$. It 
has the expression: 
$${\cal E}(\eta) \ = \ - \ \lim_{s \rightarrow 0} 
\frac{ \int_{\Omega} w(\nabla \eu(\phi_{s},\eu^{0})) \mbox{ d}x \ - \ 
\int_{\Omega} w(\nabla \eu) \mbox{ d}x}{s} \ \ .$$
With this notation the relation (\ref{apulim}) becomes: 
$$-{\cal E}(\eta) \ \leq  \ - K2(1_{\Omega},\eta,\eu^{0}) \ \ .$$
In fact we have proved more, namely that 
$$-CM^{+}(\eu)^{C}(\Omega) \ \leq \ -{\cal E}(\eta) \ \leq  \ - K2(1_{\Omega},\eta,\eu^{0}) \ \ .$$
Therefore we have the following relation: 
\begin{equation}
\mid K2 \mid (\eu(K,\eu^{0})) (\Omega) \ \leq \ \sup \left\{ {\cal E}(\eta) \mbox{ : } \eta \mbox{ satisfies (\ref{didi}) } \right\} \ 
\leq \ CM^{+}(\eu(K,\eu^{0}))^{C}(\Omega) \ \ . 
\label{minimaxi}
\end{equation}
This means that the supremum of the initial energy release rate, taken over all possible smooth propagations of the initial crack $K$ with sub-unitary 
speed, is bounded by the concentration coefficients $\mid K2 \mid$ and $CM^{+}$, 
 calculated for the displacement of the body with the initial crack $K$. 

The assumption that Griffith and Irwin criterions of brittle crack propagation are equivalent is usual in classical fracture mechanics. From our 
considerations it is straightforward that this assumption means that at any moment the coefficients $\mid K2 \mid$ and $CM^{+}$ are equal, which 
leads to the following physical statement: during the evolution, the crack takes energy from the body with maximal speed. 
We leave for further examinations whether  this statement can be 
generally proved or it has to be supposed.

\vspace{3.cm}

{\it August 1, 1997} \hfill {\it Marius Buliga} 

\hfill{\it Institute of Mathematics} 

\hfill{\it of the Romanian Academy}

\end{document}